\documentclass{article}
\usepackage{arxiv}

\usepackage[utf8]{inputenc} 
\usepackage[T1]{fontenc}    
\usepackage{hyperref}       
\usepackage{url}            
\usepackage{booktabs}       
\usepackage{amsfonts}       
\usepackage{nicefrac}       
\usepackage{microtype}      
\usepackage{lipsum}
\usepackage{graphicx}

\usepackage{url}

\usepackage{amsmath}
\usepackage{amsthm}
\usepackage{subcaption}

\title{Conics Quadrics Mapping \& Cones}
\author{Paul Zsombor-Murray \\
  McGill University \\
  Montreal QC H3A 0C3, Canada \\
  \texttt{paulzm363@gmail.com}
  \And
  Martin Pfurner \\
  University of Innsbruck\\
  6020 Innsbruck, Austria \\
  \texttt{martin.pfurner@uibk.ac.at}
}

\date{}

\begin{document}

\maketitle              

\begin{abstract}
An efficient way to get implicit equations of conics on five points
and quadrics on nine, using pencils of conics and quadrics,
is revealed. Parallel axis right cones intersect on a conic.
An example, to show how to place five coplanar points on a cone,
using kinematic mapping with dual quaternions is presented. A second
congruent cone is found as a translation of the first. Cone symmetry helps to explain how mapping produces eight real solutions, apparently
all different, belong to a unique cone pair intersection.
Future extension of this, pertaining to a nine point quadric,
can be contrived if a way to map planar points to parallel
axis cones is formulated.
\end{abstract}

\section{Conic in the Plane}\label{CiP}
Six coefficients $a_{ij}$ of the implicit equation
$a_{00}x_0^2+2a_{01}x_0x_1+2a_{02}x_0x_2+a_{11}x_1^2
+2a_{12}x_1x_2+a_{22}x_2^2=0$
of a conic on five points $PQRST$ are usually computed as $5\times 5$ sub-determinants of the matrix, Eq.~\ref{6x6}
\begin{equation}\left|\begin{array}{cccccc}
x_0^2&2x_0x_1&2x_0x_2&x_1^2&2x_1x_2&x_2^2\\
1&p_1&p_2&p_1^2&p_1p_2&p_2^2\\
1&q_1&q_2&q_1^2&q_1q_2&q_2^2\\
1&r_1&r_2&r_1^2&r_1r_2&r_2^2\\
1&s_1&s_2&s_1^2&s_1s_2&s_2^2\\
1&t_1&t_2&t_1^2&t_1t_2&t_2^2\end{array}\right|=0,\label{6x6}
\end{equation}
requiring about 3600 FLOPS. Applying LQ decomposition reduces
the number but effectively cancels the improvement. Consider an efficient alternative method using conic pencils.

\subsection{Scaled Sum of Line Equation Products}\label{1.1}
\begin{figure*}[h]
\begin{center}
\includegraphics[width=\textwidth]{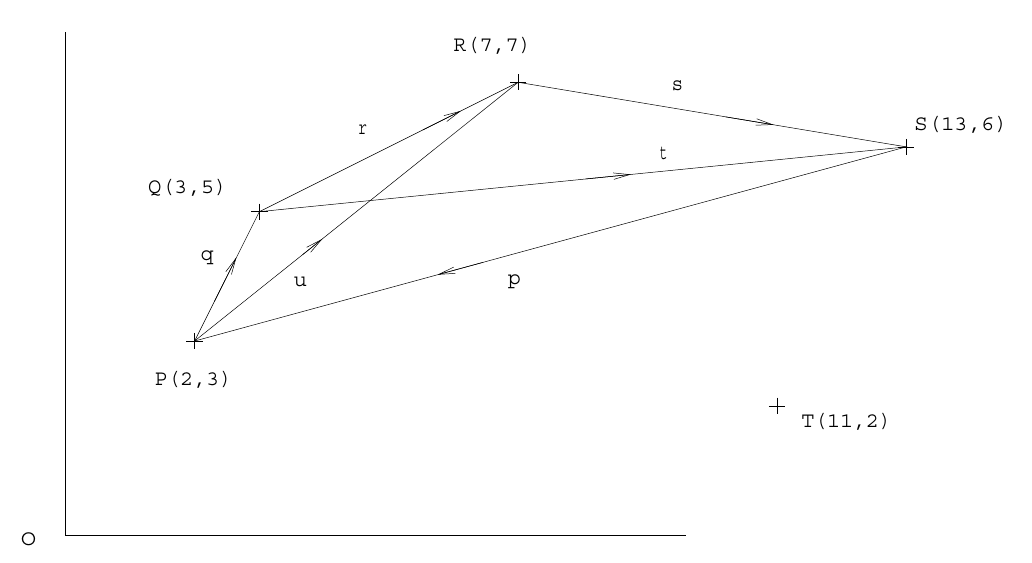}
\caption{Five Points, Four of Which Are Spanned by Six Line
Segments}
\label{5PC}
\end{center}
\end{figure*}
Fig.~\ref{5PC} shows sides and diagonals of
a quadrilateral on points $P,Q,R,S$. Lines $p,q,r,s,t,u$
are formed as
\[p=S\cap P,\;q=P\cap Q,\;r=Q\cap R,\;s=R\cap S,\;
t=Q\cap S,\;u=P\cap R.\]
The fifth point is $T$. The novel method
generates the conic equation as a scaled
$(\lambda,\mu)$ sum of binary line equation products, {\itshape viz.},
the equation Eq.~\ref{LM} built on five given points.

\begin{equation}
\lambda pr+\mu qs=0\label{LM}\end{equation}
In what follows upper case unsubscripted letters are points
with lower case subscripted homogeneous coordinates, {\itshape viz.},
$P\{p_0:p_1:p_2\}$ and lower case unsubscripted letters with
upper case subscripted homogeneous coordinates are lines,
{\itshape viz.}, $p\{P_0:P_1:P_2\}$. In some places, {\itshape e.g.},
Eq.~\ref{LM}, $p,r,q,s$ signify {\itshape line equations},
not coordinates.

Computation begins with Eqs.~\ref{XPL}. This relates to the
numerical example in Fig.~\ref{5PC} with the six lines given by
their implicit equations. Each of 12 coefficients takes
three FLOPs to calculate for a total of 36 thus far.
\begin{equation}
p\{P_0:P_1:P_2\}=\left[\begin{array}{c}1\\s_1\\s_2\end{array}
\right]{\mathbf\times}\left[\begin{array}{c}1\\p_1\\p_2\end{array}
\right]\equiv\left[\begin{array}{c}1\\13\\6\end{array}\right]
{\mathbf\times}\left[\begin{array}{c}1\\2\\3\end{array}\right]
\rightarrow 27+3x_1-11x_2=0\label{XPL}\end{equation}
Similarly
\[\begin{array}{c}
q\rightarrow 1-2x_1+x_2=0,\;\;r\rightarrow -14-2x_1+4x_2=0,\;\;
s\rightarrow -49+x_1-6x_2=0\\
t\rightarrow -47-x_1+10x_2=0,\;\;u\rightarrow-7-4x_1+5x_2=0.
\end{array}\]
Line pair $t,u$ is included to show the line equation products/sum equivalence $pr+qs=tu$. {\itshape I.e.}, using the planes with numerical
coefficients in Eq.~\ref{XPL}, obtained with coordinates
from Fig.~\ref{5PC},
\[\begin{array}{c}
(27+3x_1-11x_2)(-14-2x_1+4x_2)-(1-2x_1-x_2)(-49+x_1+6x_2)\\
-(-47-x_1+10x_2)(-7-4x_1+5x_2)=0. \end{array}\]
To prove this equivalence, some liberty was taken with
$\pm$ signs due to homogeneity of plane coordinate and
non-commutativity of cross products. Notice how the arrows
$p\rightarrow q\rightarrow r\rightarrow s$ form a closed clockwise
sequence. Including the pair $ts$ with any other, say $pr$,
breaks it. This is important. It helps
to explain univariate polynomial multiplicity.
This is encountered later, Eq.~\ref{3UVP},
in mapping five planar points to a right cone to show how this
ties into the elementary concept of conic section.
Let us now turn to finding $\lambda,\mu$ so that the
un-factorable conic {\itshape does} contain $T$.

 \subsection{Multipliers \& Conic}
\begin{figure*}[h]
\begin{center}
\includegraphics[width=.5\textwidth]{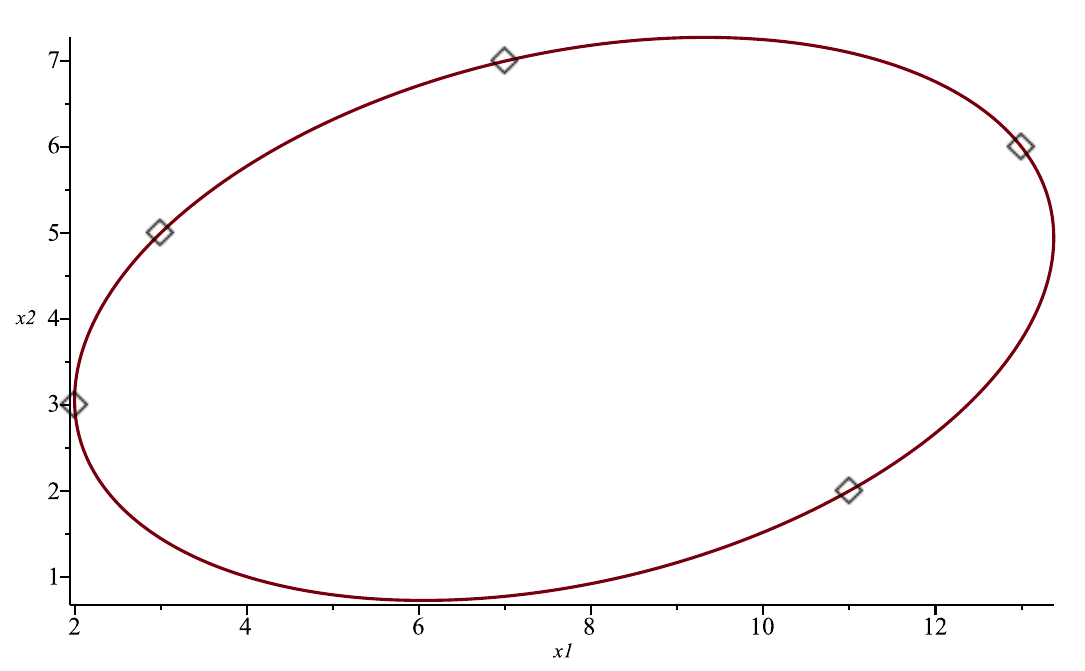}
\caption{Verification of Eq.~\ref{Ceq}}
\label{5Pcn52c}
\end{center}
\end{figure*}

Reproducing $\lambda pr +\mu qs=0$ with numerical values from
Fig.~\ref{5PC} and using Eqs.~\ref{XPL} substituted for $X\{x_0:x_1:x_2\}$ gives us Eq.~\ref{subT} for $\lambda,\mu$.
\begin{equation}\begin{array}{c}
\lambda(P_0x_0+P_1x_1+P_2x_2)(R_0x_0+R_1x_1+R_2x_2)\\
+\mu(Q_0x_0+Q_1x_1+Q_2x_2)(S_0x_0+S_1x_1+S_2)=0\rightarrow\\
\lambda(27+3t_1-11t_2)(-14-2t_1+4t_2)+\mu(1-2t_1+t_2)(-49+t_1+6t_2)
=0\rightarrow\medskip\\
\lambda(27+33-22)(-14-22+8)+\mu(1-22+2)(-49+11+12)\\=
-1064\lambda+494\mu\rightarrow\lambda=494,\;\;\mu=1064
\end{array}\label{subT}\end{equation}
This accounts for about 22 or 23 more FLOPs, about 2\% of that
required to expand six $5{\mathbf\times}5$ determinants. The
conic equation is Eq.~\ref{Ceq} while Fig.~\ref{5Pcn52c} is
that conic.
\begin{equation}
-238868x_0^2+57912x_0x_1+83676x_0x_2-5092x_1^2+5092x_1x_2-15352x_2^2
=0\label{Ceq}\end{equation}

\section{Nine Point Quadric}
\label{nine-point-quadric}
Replacing Eq.~\ref{6x6} with a $10{\mathbf\times}10$
matrix and evaluating ten $9{\mathbf\times}9$ determinants
incurs at least an order of magnitude greater number
of FLOPs than the conic. Now quadric pencils replace conic
pencils. An octahedron replaces
the quadrilateral in Fig.~\ref{5PC}.
Six points on its vertices and eight planes, each defined by
three of the points, together with three remaining points
together form a system of four plane pairs. Two plane pairs
provide a double covering of the vertices. Combining the
remaining points requires four multipliers
$\alpha,\beta,\gamma,\delta$ as opposed to the two,
$\lambda,\mu$, for conics.

\begin{itemize}
\item Summarizing, three octahedron vertices among
$ABCDEF$, {e.g.}, $ABC$, define a plane,
\item And, automatically, its partner that effects
complete covering. Upper/lower case and Greek letters
apply to point/plane in the same convention as point/line
was applied to conics in Section~\ref{1.1}.
\item Four of these are chosen from the available 10
plane pairs shown below. The choice is $_{10}C_4=210$.
\item From these
\[\begin{array}{c}
ABC\leftrightarrow DEF,\;\;ABD\leftrightarrow CEF,\;\; ABE\leftrightarrow CDF,\;\;ABF\leftrightarrow CDE,\\ ACD\leftrightarrow BEF,\;\;ACE\leftrightarrow BDF,\;\; ACF\leftrightarrow BDE,\;\;ADE\leftrightarrow BCF,\\ ADF\leftrightarrow BCE,\;\;AEF\leftrightarrow BCD\end{array}\]
\item Four plane pairs
\[\begin{array}{c}ABF\leftrightarrow CDE\equiv pv,\;\;
ADE\leftrightarrow BCF\equiv ru,\\
ADF\leftrightarrow BCE\equiv ts,\;\;
ABE\leftrightarrow CDF\equiv qw\end{array}\]
are chosen.
\item Like the two of three line pairs in Eq.~\ref{LM}
the four plane pairs, arbitrarily selected, are assembled
into the sum of four plane products and multipliers, Eq.~\ref{4e},
\begin{equation}
pv\alpha+ru\beta+ts\gamma+qw\delta=0.\label{4e}\end{equation}
\item  Any other choice of four plane pairs, {\itshape e.g.},
\[
ABF\leftrightarrow CDF,\;\;AED\leftrightarrow EBC,\;\;
ADF\leftrightarrow BCF,\;\;ECD\leftrightarrow AEB,\]
substituted into Eq.~\ref{4e} yields the same quadric.
\end{itemize}

\subsection{Eight Chosen Plane Pairs}
\begin{figure*}[h]
\begin{center}
\includegraphics[width=\textwidth]{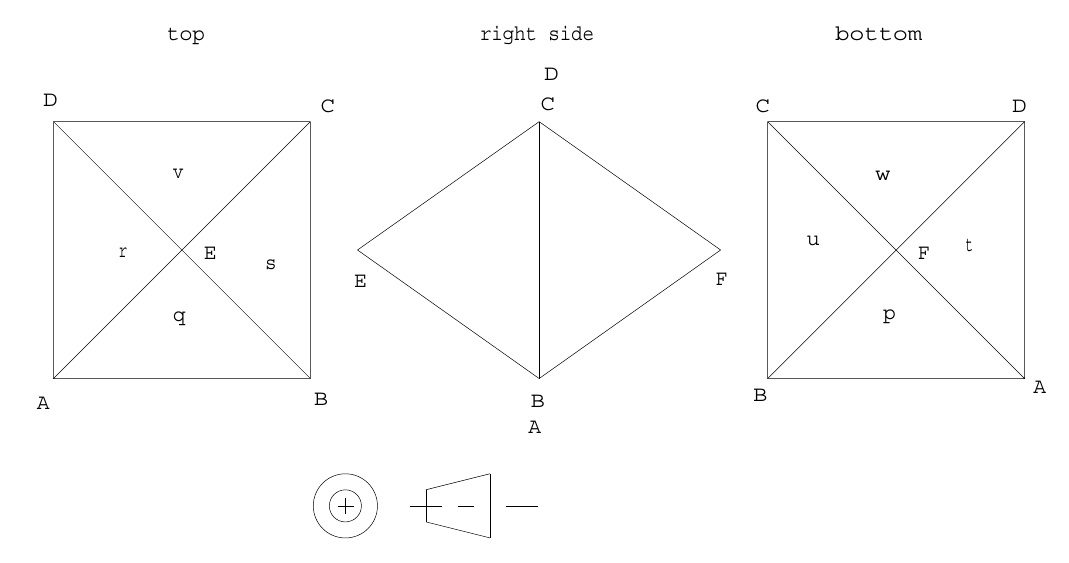}
\caption{Point \& Plane Topology of a Regular Octahedron}
\label{TOPO52D}
\end{center}
\end{figure*}

Rather than trying to show the octahedron, like the
actual quadrilateral with five given points from which
two line pairs were chosen from Fig.~\ref{5PC},
the nine points for a numerical example
\[\begin{array}{c}
A(4,3,0),\;\;B(4,10,4),\;\;C(9,3,8),\;\;D(-1,7,2),\;\;
E(2,3,5),\;\;F(-3,9,7),\\
G(12,0,0),\;\;H(0,10,0),\;\;J(0,0,5)\end{array}\]
are referred to Fig.~\ref{TOPO52D}, a regular octahedron.
This is deemed better than trying to convey the relation
with six points on eight scalene triangular faces, some
possibly internal.

\subsection{Numerical Example of Nine Point Quadric}
\begin{figure*}[h]
\begin{center}
\includegraphics[width=.5\textwidth]{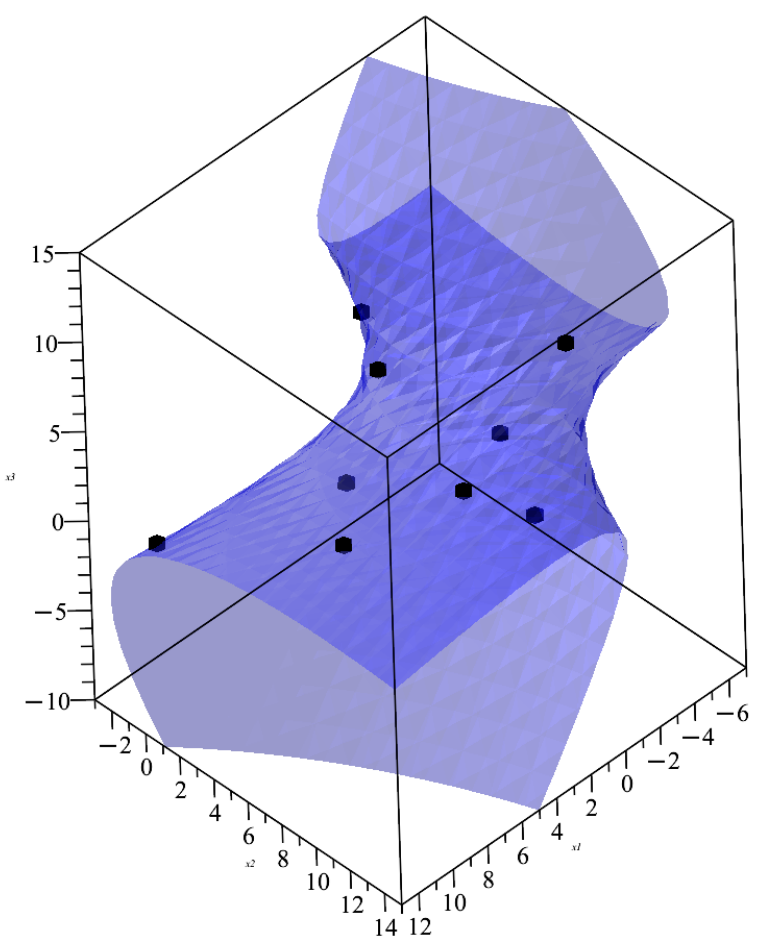}
\caption{Hyperboloid of One Sheet Defined by
the Nine Given Points}
\label{BluH1S}
\end{center}
\end{figure*}

Given our choice of four plane
pairs the quadric equation with multipliers looks like
Eq.~\ref{Qsym}.
\begin{equation}\begin{array}{c}
\alpha(P_0x_0+P_1x_1+P_2x_2+P_3x_3)(V_0x_0+V_1x_1+V_2x_2+V_3x_3)
\medskip\\
+\beta(R_0x_0+R_1x_1+R_2x_2+R_3x_3)(U_0x_0+U_1x_1+U_2x_2+U_3x_3)
\medskip\\
+\gamma(T_0x_0+T_1x_1+T_2x_2+T_3x_3)(S_0x_0+S_1x_1+S_2x_2+S_3x_3)
\medskip\\
+\delta(Q_0x_0+Q_1x_1+Q_2x_2+Q_3x_3)(W_0x_0+W_1x_1+Q_2x_2+Q_3x_3)
=0\end{array}\label{Qsym}\end{equation}
We get the three necessary equations, Eqs.~\ref{3neq},
by computing
plane coordinates from the points and successively substituting
\[
G\{g_0:g_1:g_2:g_3\},H\{h_0:h_1:h_2:h_3\},J\{j_0:j_1:j_2:j_3\}\]
for point $X\{x_0:x_1:x_2:x_3\}$
in Eqs.~\ref{Qsym} to yield three equations, Eqs.~\ref{3neq}.
\begin{equation}\begin{array}{c}
-103944\alpha+9700\beta+7410\gamma-109136\delta=0\\
3552\alpha+6968\beta+15936\gamma-38612\delta=0\\
-15343\alpha-9167\beta-3699\gamma+19918\delta=0
\end{array}\label{3neq}\end{equation}
These are solved with the Grassmannian determinant
in Eq.~\ref{GmDt} to yield
$\alpha,\beta,\gamma,\delta$.
\begin{equation}\left|\begin{array}{cccc}
\alpha&\beta&\gamma&\delta\\-103944&9700&7410&-109136\\
3552&6968&15936&-38612\\-15343&-9167&-3699&19918\end{array}
\right|\rightarrow
\begin{array}{c}
\alpha=-9842336242680\\
\beta=39532196597640\\
\gamma=19311493179280\\
\delta=14198910257520\end{array}\label{GmDt}\end{equation}
These in turn go back into Eq.~\ref{Qsym} as Eq.~\ref{Qeq}
and plot as Fig.~\ref{BluH1S} together with the nine points thereon.
\begin{equation}\begin{array}{c}
-202095981901413x_1^2+22422685213194x_1x_2-1191822696049068x_1x_3\\
-158099339159010x_2^2-558476852988570x_2x_3-740341605937509x_3^2\\
+4710658547758491x_1+4323601509519942x_2+9186924265547229x_3\\
-27426081179298420=0\end{array}\label{Qeq}\end{equation}

\section{Five Points Placed on a Conic Section}
Given five coplanar points $ABCDE$ --different point names for a
different numerical example to exploit advantage of their ideal
placement without too much sacrifice of generality--
that define a unique conic according to Section~\ref{nine-point-quadric}, above,
what is the relation of $\mu :\lambda$ in
Eq.~\ref{LM} to a planar section of an right cone on the five
points, chosen to simplify algebra, with apex on origin,
angle $\pi/2$ and axis on the $z$-axis? A possible approach,
but as yet without answering the question, might be
to find the disposition of that plane and the five points by
mapping them to an actual cone using dual quaternion (DQ)
algebra as outlined below.
\begin{itemize}
\item Define a planar rigid body in the end effector frame EE
containing\newline $A(0,0,0),B(b_1,0,0),C(c_1,c_2,0),D(d_1,d_2,0),E(e_1,e_2,0).$
\item Define the fixed frame FF with cone apex and axis as
specified above.
\item Map the five points in EE to FF as follows.
\item Place $A$ on the cone generator defined by the intersection
of planes $y=0$ and $x-z=0$.
\item Map the remaining four points to the cone
$x^2+y^2-z^2=0$ in FF with a triple DQ product.
\item Placement of $A$ provides two linear constraint equations.
\item Mapping of $BCDE$ entails four quadratic constraints.
\item Two more constraints are obtained with norming
and Study conditions.
\[x_0^2+x_1^2+x_2^2+x_3^2-1=0,\;\;x_0y_0+x_1y_1+x_2y_2+x_3y_3=0\]
\item These eight constraints are sufficient to extract the
dual quaternion necessary to place the five point rigid body array
on the cone, more that enough to define the pose of the planar section so it can be compared with the $\mu :\lambda$ ratio.
\end{itemize}

\subsection{Dual Quaternions \& Multiplications}
DQ ${\mathbf Q}$ pre-multiplies typical point
quaternion
${\mathbf C}$ in EE post-multiplied by DQ ${\mathbf Q}'_\epsilon$
to produce image of ${\mathbf C}'$ in FF, Eq.~\ref{DQM},
using the procedure detailed in \cite{GF22}. ${\mathbf Q}'_\epsilon$
is conjugate of ${\mathbf Q}_\epsilon$ which is ${\mathbf Q}$ with
all $y_i$ negated, $-y_i$.

{\small \begin{equation}\begin{array}{c}
\left[\begin{array}{c}x_0\\x_1\\x_2\\x_3\\y_0\\y_1\\y_2\\y_3
\end{array}\right]\left[\begin{array}{c}1\\0\\0\\0\\0\\c_1\\c_2\\0
\end{array}\right]\left[\begin{array}{c}
x_0\\-x_1\\-x_2\\-x_3\\-y_0\\y_1\\y_2\\y_3\end{array}\right]
=\left[\begin{array}{c}1\\0\\0\\0\\0\\c'_1\\c'_2\\c'_3
\end{array}\right]
=2\left(\left[\begin{array}{c}0\\0\\0\\0\\0\\
c_1(x_0^2+x_1^2)-c_2(x_0x_3-x_1x_2)\\
c_1(x_0x_3+x_1x_3)+c_2(x_0^2+x_2^2)\\
-c_1(x_0x_2-x_3x_1)+c_2(x_0x_1+x_2x_3)\end{array}\right]\right.
\medskip\\
+\left.\left[\begin{array}{c}0\\0\\0\\0\\0\\
x_0y_1-x_1y_0+x_2x_3-x_3x_2\\x_0y_1-x_2y_0+x_3y_1-x_1y_3\\
x_0y_3-x_3y_0+x_1y_2-x_2y_1\end{array}\right]\right)
+\left[\begin{array}{c}1\\0\\0\\0\\0\\-c_1\\-c_2\\0\end{array}
\right]\equiv{\mathbf Q}{\mathbf C}{\mathbf Q}'_\epsilon ={\mathbf C}'\end{array}\label{DQM}\end{equation}}
It is noted with satisfaction that the result is consistent with that
of a point represented as dual quaternion.
$C(c_1,c_2,0)\rightarrow C'(c'_1,c'_2,c'_3)$
describes the desired mapping of typical points $C,C'$ expressed in
Cartesian coordinates.

\subsection{Numerical Examples}\label{NX}
\begin{figure}[h]
  \centering
  \begin{subfigure}{0.25\textwidth}
    \includegraphics[width=\textwidth]{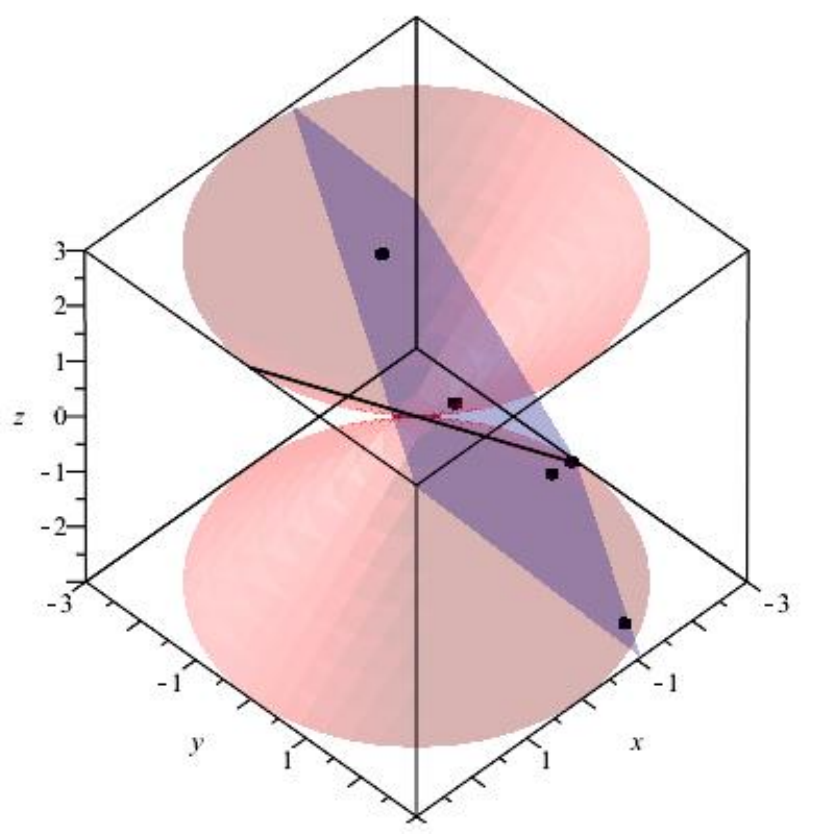}
    \caption{$x_0=0.1380$, the Cone $x^2+y^2-z^2=0$, 5 Points}
  \end{subfigure}
  \hfill
  \begin{subfigure}{.4\textwidth}
    \includegraphics[width=\textwidth]{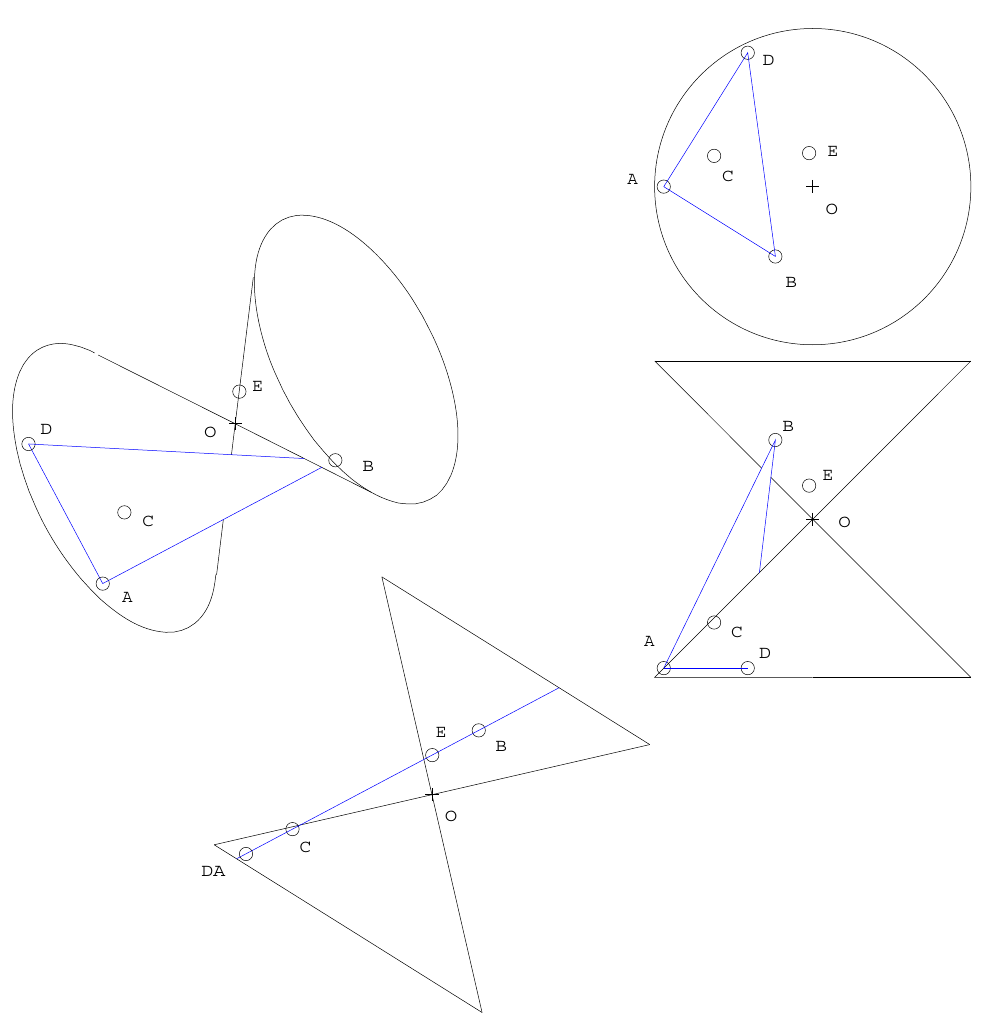}
    \caption{Three Views of Cone \& 3-Point Plane Segment}
  \end{subfigure}
  \hfill
  \begin{subfigure}{.25\textwidth}
    \includegraphics[width=\textwidth]{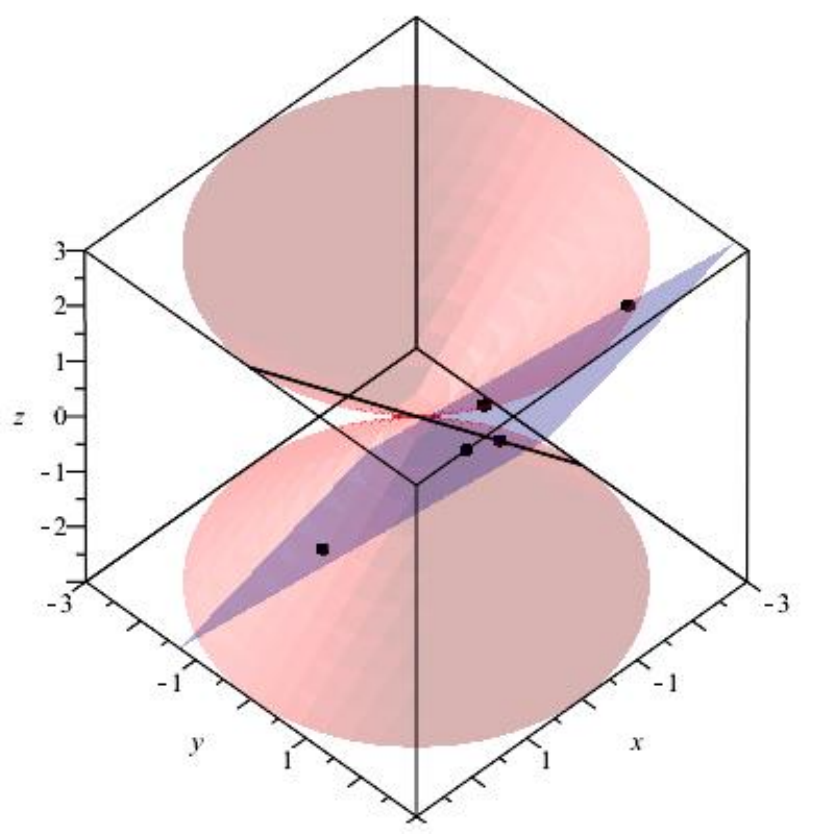}
    \caption{$x_0=0.6333$, the Cone $x^2+y^2-z^2=0$, 5 Points}
  \end{subfigure}
\caption{Cones Points \& Planes}
\label{3CnFgs}
\end{figure}

The five points $ABCDE$ in EE are given in homogeneous coordinates as
\[\begin{array}{c}
A\{1:0:0:0\},B\{1:5:0:0\},C\{1:1:-1:0\},\\D\{1:0:-3:0\},E\{1:4:-2:0\},
\end{array}\]
mapped to dual quaternions and pre- and post-multiplied as shown in
Eq.~\ref{DQM} to yield $A',B',C',D',E'$ in FF in terms of
eight Study parameters $x_i,y_i$. Two equations are provided by the
simplifying condition of placing $A$ on the generator of a right cone,
apex on the chosen origin in FF, axis on FF $z$-axis and apex angle of
$\pi/2$, on the intersection of planes $y=0$ and the tangent plane
with homogeneous coordinates $\{0:1:-1:0\}$. Four more arise
with the transformation of $B,C,D,E$ according to Eq.~\ref{DQM}.
Adding Study quadric and norming condition make
eight. These are readily solved with a Gr\"obner basis
algorithm like that found in {\itshape Maple} software. In this case
{\ttfamily plex}-order was chosen to yield a univariate polynomial (UVP)
in $x_0$ as the first basis followed by seven bivariate bases,
all seven linear in $x_i,y_i$ except $x_0$ available from UVP.
For this example the UVP is of degree 32 but
factors in three of 8, 8 and 16 as shown in Eqs.~\ref{3UVP}.
\begin{equation}\begin{array}{l}
(1183744x_0^8-591872x_0^6+97600x_0^4-5904x_0^2+81=0)\\
(111183744x_0^8-1775616x_0^6+878400x_0^4-159408x_0^2+6561=0)\\
(516716384256x_0^{16}-10233432768512x_0^{14}\\+7827661127680x_0^{12}
-1021055926272x_0^{10}\\
-834991220736x_0^8+75713882112x_0^6+209299178880x_0^4\\
+60512832x_0^2+6561=0)\end{array}\label{3UVP}\end{equation}
The factor of degree 16 contains only
complex roots in $\pm$ pairs. A DQ vector is homogeneous.
Negation does not alter its effect. The two factors of degree 8 also
come in $\pm$ pairs so each accounts for only 4 real solutions.
Approximate solutions for $x_0$ from the two factors of degree 4
are tabulated below.
\[\begin{array}{rrrr}
+\mathbf{0.1380}&\pm 0.3411&\pm 0.3656&\pm 0.4806\medskip\\
\pm 0.2391&\pm 0.5907&+\mathbf{0.6333}&\pm 0.8324\end{array}\]
Using the two boldfaced positive values of $x_0$ in the preceding
table the two resulting DQ are shown below. Using the first,
$x_0=0.1380$, and the DQ derived therefrom
\[\left[0.1389\;\;0.08324\;\;-0.2391\;\;0.4806\;\;1.8555\;\;-0.5330\;\;
-0.4972\;\;0.1428\right]^\top\]
the five transformed points are
\[\begin{array}{c}
A'(-2.8265,0,-2.8265),B'(-0.7076,-1.3267,1.5036),\\
C'(-1.8720,0.5822,-1.9605),
D'(-1.2344,2.5427,-2.8265),\\E'(-0.0700,0.6337,0.6376).
\end{array}\]
These results are shown plotted on a right cone and the plane
 on which they lie in Fig.~\ref{3CnFgs} (a).
Homogeneity and only even powers of $x_0$ appearing
in UVPs explains how
$\pm$ solution pairs reduces number of real ones in
both factors from 8 to 4. Symmetry, indicating solutions
above and below cone vertex on the same generator, shows
that the 4 are really 2.
What about the second, similar, factor of degree 8 that
can also be argued down from 8 to 2 solutions? These are
rigidly paired with the first 8, shared by a congruent
cone displaced from the first by
translation along the generator which maintains $A'$ on it.
The translation, described in~\ref{sc3.3}, of the
first cone produces the second which intersects it to share
a plane, the 5 points and the conic because it is known
\cite{PA18} that {\itshape any} conic section can be generated
as the intersection of parallel axis right cones.
Finally to explain remaining 2 real solutions note that
coordinates of points $ABCDE$ in EE are on plane $z=0$.
In the plane, the face triangles $ABC\neq ACB$ are not
congruent. In space $ABC=ACB$ they can be made to transform
as such \footnote{This is indeed the situation in the case of both
DQ transformations, {\itshape i.e}, that illustrated by
Figs~\ref{3CnFgs} (a), (b) and (c).
Section~\ref{MRG} (Appendix) explains how mirroring happens.}.

It is the ultimate aim of this work to relate the ratio $\mu : \lambda$, encountered previously,
to the displacement between parallel axis right cones vertices and,
if possible, to extend it to four parameter, $\alpha,\beta,\gamma,\delta$, situation case of quadrics. Both
are topics for future investigation.

For now look at the results obtained with $x_0=0.6333$. The DQ is
\[
\left[0.6333\;\;0.3411\;\;-0.3656\;\;0.5907\;\;0.7005\;\;
-0.7510\;\;0.1877\;\;-0.2012\right]^\top . \]
The five transformed points are
\[\begin{array}{c}
A'(-1.5036,0,-1.5036),B'(-1.3301,2.4940,2.8265),\\
C'(-0.4713,0.4294,-0.6376),D'(1.4891,-0.2082,-1.5036),\\
E'(0.6304,1.8564,1.9605). \end{array}\]
These are plotted with cone and plane in Fig.~\ref{3CnFgs} (c).
This is an image that shows the translated cone. The five points
$A'B'C'D'E'$ are on it and the plane intersecting both cones.

\subsection{Right Cones with Parallel Axes}
Recall Eq.~\ref{3UVP} and the 32 solutions for $x_0$ it produces.
The first two factors yield 16 real values for $x_0$.
The third factor has 16 complex roots. They are artifacts of absolute points on intersecting coplanar circles ruling two congruent cones
with parallel axes, the circles being normal to cone axes.
Difference of circle equations is that of the real line between
the two. It is the purpose here to show how the second, displaced,
cone intersects the origin centred one on a shared conic. Any one
of the real solutions yields a cone that can then be translated to
get the second. {\itshape I.e.}, taking Fig.~\ref{3CnFgs} (a) to be the
origin centred one and Fig.~\ref{3CnFgs} (c) the displaced one the
plane (blue) in each image is the same, but on the second cone
it has been translated downward and rotated. Fig.~\ref{3CnFgs} (b)
shows projections of Fig.~\ref{3CnFgs} (a) with plane segment $ABC$
thereon. Fig.~\ref{3CnFgs} (a) and (b) were produced by the
mapping that used $x_0=0.1380$ while Fig.~\ref{3CnFgs} (c) used
$x_0=0.6333$.

\subsection{Cones \& Planes}\label{sc3.3}
\begin{figure*}[h]
\begin{center}
\includegraphics[width=.5\textwidth]{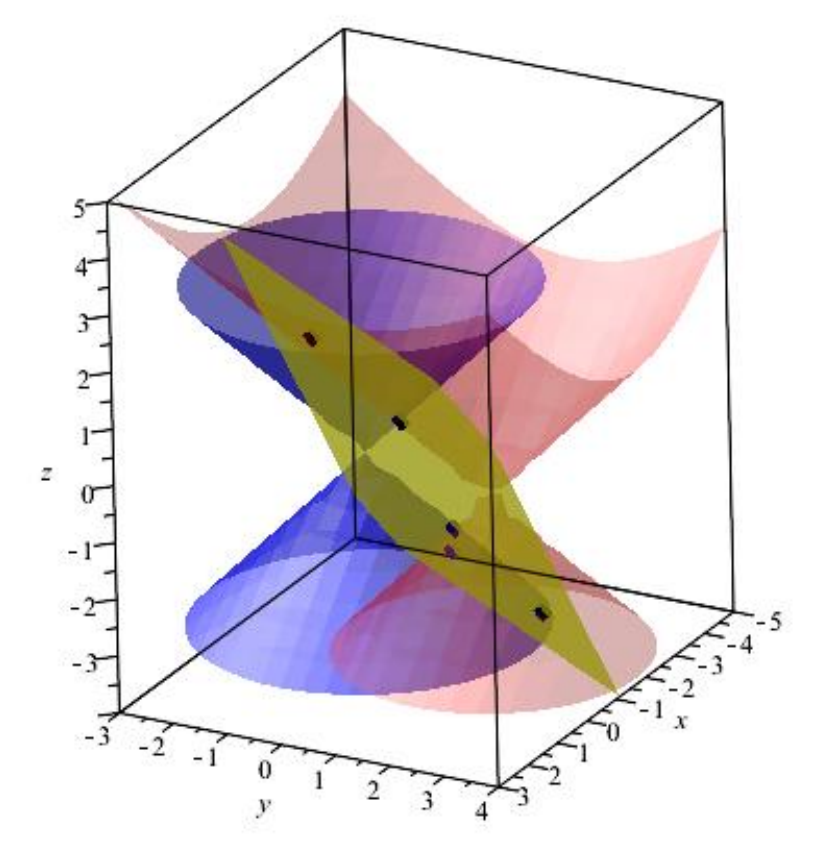}
\caption{Five Points \& Four Surfaces}
\label{3F5P61u}
\end{center}
\end{figure*}

The translated cone equation is given by Eq.~\ref{TK}.
\begin{equation}\begin{array}{c}
\left[1\;\;x\;\;y\;\;z\right]\left[\begin{array}{cccc}
1&-t_1&-t_2&t_3\\0&1&0&0\\
0&0&1&0\\0&0&0&1\end{array}\right]
\left[\begin{array}{cccc}0&0&0&0\\0&1&0&0\\0&0&1&0\\0&0&0&-1
\end{array}\right]
\left[\begin{array}{cccc}1&0&0&0\\-t_1&1&0&0\\-t_2&0&1&0\\
-t_3&0&0&1\end{array}\right]
\left[\begin{array}{c}1\\x\\y\\z\end{array}\right]=\medskip\\
(t1^2+t_2^2-t_3^2)-2t_1x-2t_2y+2t_3z+x^2+y^2-z^2=0
\end{array}\label{TK}\end{equation}
Subtracting this from the origin centred cone $x^2+y^2-z^2=0$
gives Eq.~\ref{KX}, the plane of intersection.
\begin{equation}
t_1^2+t_2^2-t_3^2-2t_1x-2t_2y+2t_3z=0\label{KX}
\end{equation}
Substituting
\[\begin{array}{c}A'(-2.8265,0,-2.8265),B'(-0.7076,-1.3267,1.5036),\\
C'(-1.8720,0.5822,-1.9605)\end{array}\]
into Eq.~\ref{KX} subtracting the first from the other two gives
a quadratic in $t_1,t_3$ and two linear equations in $t_1,t_2,t_3$,
Eqs.~\ref{3t}.
\begin{equation}\begin{array}{c}
t_1^2+t_2^2-t_3^2+5.6530t_1-5.6530t_3=0\\
4.2378t_1-2.6534t_2-8.6602t_3=0\\
1.9090t_1+1.1644t_2-1.7320t_3=0\end{array}\label{3t}
\end{equation}
Solution gives the origin centred cone with $t_1=t_2=t_3=0$ and
the displaced cone with $t_1=-1.9418,t_2=1.2160,t_3=-1.3227$
whose equation is Eq.~\ref{eqk2}. The plane of intersection is
the same with squared terms removed.
\begin{equation}
3.5+3.884x-2.432y-2.6458z+x^2+y^2-z^2=0\label{eqk2}\end{equation}
Using $A',B',C'$ from above, plane equation obtained with mapped
points is Eq.~\ref{km}.
\begin{equation}3.3069+3.6699x-2.2981y-2.5000z=0\label{km}
\end{equation}
The two cones, two planes --which are the same
\footnote{Blue and yellow overlap to produce what appears
to be a single green plane.}-- and the five
given points appear in Fig.~\ref{3F5P61u}.

\section{Conclusion}
\begin{itemize}
\item Efficient ways to compute implicit equations of conics
and quadrics on points and the underlying geometry were described.
\item Kinematic mapping was used to move planar point arrays onto
conic sections, relating these to their equations.
\item How the section is shared on the intersection curve between
two congruent cones was shown.
\item Multipliers extract the unique conic or quadric using two
singular members of the pencil in which it resides.
\item Future work is intended to connect multipliers with the
displacement between intersecting cone pair.
\item To do this for quadrics can {\itshape Zyklographie} of planar
points and cones be extended to comprehend spatial points
and four dimensional cones?
\end{itemize}

\section{Bibliography}
The note, \cite{GF22}, is useful
in navigating the intricacies of dual
quaternion algebra. \cite{PA18} deals with
conics on quadrics. Pencils and other aspects of conics and
quadrics are dealt with in \cite{UC16},\cite{UQ20},\cite{XGM27},\cite{CGM26},\cite{PPG99}.
\cite{SAG84} relates parametric and implicit equations
to elements of simultaneous polynomials.

\newpage

\section{Appendix}\label{MRG}
\begin{figure*}[h]
\begin{center}
\includegraphics[width=.7\textwidth]{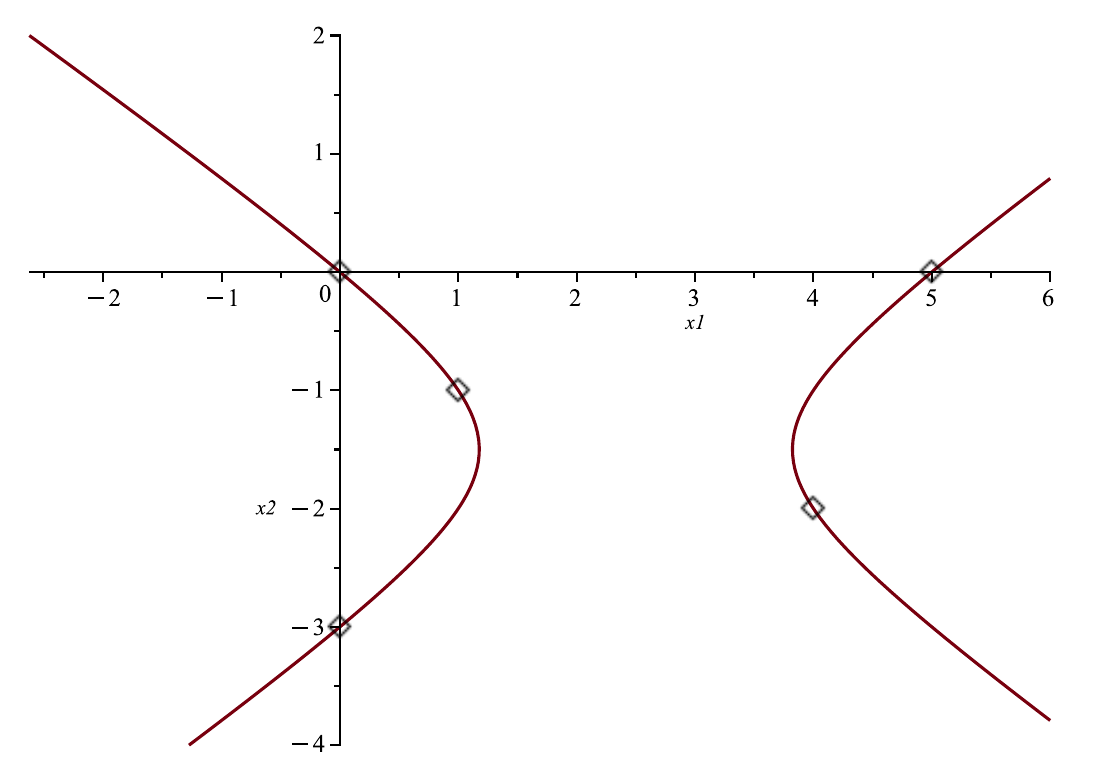}
\caption{Hyperbola on Five Points}
\label{Hy}
\end{center}
\end{figure*}

Spurious solutions are generated because the five coplanar points,
presented in Section~\ref{NX} as homogeneous spatial coordinates
on plane $z=0$, are also used to find the hyperbola they represent
using the pencil of conics method described in Section~\ref{CiP} where points
are used in cross-products as homogeneous {\itshape planar} coordinates.
Because of the sequence in which points on the quadrilateral $ABCD$
are presented the result shown in Fig.~\ref{Hy} has $AB$ from origin along positive $x$-axis while $AD$ extends from origin along
negative $y$-axis. In the second auxiliary view in Fig,~\ref{3CnFgs}
(b) $AD$ is along the positive $y$-axis. DQ transformation has chosen to legitimately induce a half-turn  about an axis in plane $z=0$
resulting in extra solutions.

\end{document}